\input phyzzx
\scrollmode
\catcode`\@=11
\def\PH@SR@V{\doubl@true \baselineskip=18.1pt plus 0.1pt minus 0.1pt
             \parskip= 3pt plus 2pt minus 1pt }
\def\PHYSREV{\papers\PhysRevtrue\PH@SR@V}
\let\physrev=\PHYSREV
\def\NAME{\author{J. A. de Azc\'arraga and J. C. P\'erez Bueno}}
\def\FTUV{\address{Departamento de F\'{\i}sica Te\'{o}rica and IFIC,\break
      Centro Mixto Univ. de Valencia-CSIC\break 
	46100-Burjassot (Valencia), Spain}}
\font\black=msbm10 scaled\magstep1
\def\bumpfootsymbolcount{\rel@x
   \iffrontpage \bumpfootsymbolpos \else \advance\lastf@@t by 1
     \ifPhysRev \bumpfootsymbolpos \else \bumpfootsymbolpos \fi \fi
   \gl@bal\lastf@@t=\pagen@ }
\def\title#1{\vskip\frontpageskip \titlestyle{\seventeenbf #1}
 \vskip\headskip }
\def\abstract{\par\dimen@=\prevdepth \hrule height\z@ \prevdepth=\dimen@
   \vskip\frontpageskip\centerline{\twelverm ABSTRACT}\vskip\headskip }
\def\ack{\par
      \ifnum\the\lastpenalty=30000\else \penalty-100\smallskip \fi
      \noindent{\bf Acknowledgements}:\quad\enspace \vadjust{\penalty5000}}
\def \dfxi {({\Number {1}}.4)}
\def\JMP#1(#2){\journal J. Math. Phys. &#1(#2)}
\def\PHL#1(#2){\journal Phys. Lett. &#1(#2)}
\def\JA#1(#2){\journal J. Algebra &#1(#2)}
\def\LMP#1(#2){\journal Lett. Math. Phys. &#1(#2)}
\def\JPH#1(#2){\journal J. Phys. &#1(#2)}
\def\CMP#1(#2){\journal Commun. Math. Phys. &#1(#2)}

\REF\APB{J. A. de Azc\'arraga and J. C. P\'erez Bueno \JMP 36 (95) 6879}

\REF\LNRT {J. Lukierski, A. Nowicki, H. Ruegg and V.N. Tolstoy \PHL B264 (91) 
331; J. Lukierski, H. Ruegg, and V.N. Tolstoy, 
{\it $\kappa$-quantum Poincar\'{e} 1994}, 
in {\it Quantum groups: formalism and applications}, 
J. Lukierski, Z. Popowicz and J. Sobczyk eds., PWN (1994), p. 359}

\REF\AACON{V. Aldaya and J.A. de Azc\'{a}rraga, {\it Int. J. of Theor. Phys.}
{\bf 24} (1985), 141}

\REF\AI{J. A. de Azc\'{a}rraga and J. M. Izquierdo, 
{\it Lie algebras, Lie groups cohomology and some applications in physics}, 
Camb. Univ. Press (1995)}

\REF\APBFEST{J. A. de Azc\'arraga and J. C. P\'erez Bueno, {\it Contractions, 
Hopf algebra extensions and covariant differential calculus}, 
FTUV 96--2 / IFIC 96--3 (q-alg/9602016),
to appear in J. Lukierski Festschrift, B. Jancewicz and J. Sobczyk eds.}

\REF\GLODZ{S. Giller, P. Kosi\'nski, M. Majewski, P. Ma\'slanka and J. Kunz 
\PHL B286 (92) 57}

\REF\MB{S. Majid \JA 130 (90) 17; {\it Israel J. Math.} {\bf 72} (1990), 133}

\REF\Praga{J. A. de Azc\'arraga and J. C. P\'erez Bueno 
{\it Czech. J. Phy.} {\bf 46} (1996), 137}

\REF\SIT{A. Sitarz \PHL B349 (95) 42}

\REF\MR{S. Majid and H. Ruegg \PHL B334 (94) 348} 

\REF\LUKRUE{J. Lukierski and H. Ruegg \PHL B329 (94) 189}

\REF\SZ{S. Zakrzewski \JPH A27 (94) 2075}

\REF\MAJSOB{S. Majid and Ya. S. Soibelman \JA 163 (94) 68; 
{\it Czech. J. Phys.} {\bf 44} (1994), 1059}

\REF\DOI{Y. Doi, {\it Commun. in Alg.} {\bf 17}(12) (1989), 3053}

\REF\LL{V. Bargmann, {\it Ann. Phys.} {\bf 59} (1954), 1; 
see also, 
J. M. L\'evy-Leblond, {\it Galilei group and Galilean invariance}, in 
{\it Group theory and its applications} vol. II, E. M. Loebl ed., Acad. Press 
(1971), p. 221}

\REF\ZAUGG{P. Zaugg, {\it The quantum Poincar\'e group from quantum 
contraction}, MIT-CPT-2353 (hep-th/9409100)}

\REF\GGKM{S. Giller, G. Gonera, P. Kosi\'nski and P. Ma\'slanska, {\it The 
quantum Galilei group} (q-alg/9505007)}

\REF\Ohn{C. Ohn \LMP 30 (94) 87}

\REF\AKS{A. Aghamohammadi, M. Korrami and A. Shariati \JPH A28 (95) L225}

\REF\KSAA{M. Korrami, A. Shariati, M. Abolhasani and A. Aghamohammadi,
{\it Mod. Phys. Lett.} {\bf A10} (1995), 873}

\REF\BCHOS{A. Ballesteros, E. Celeghini, F. J. Herranz, M. A. del Olmo and
M. Santander \JPH A28 (95) 3129}

\REF\BHOPS{A. Ballesteros, F. J. Herranz, M. A. del Olmo, C. M. Pere\~na
and M. Santander \JPH A28 (95) 7113}

\REF\CGSTa{E. Celeghini, R. Giachetti, E. Sorace, and M. Tarlini, {\it 
Contractions of quantum groups}, in Lec. Notes Math. {\bf 1510}, (1992) p. 
221; {\it J. Math. Phys.} {\bf 32} (1991), 1155}

\REF\FG{C. Fronsdal and A. Galindo \LMP 27 (93) 29;
{\it The universal $T$-matrix}, UCLA-93-TEP-2 (1993)}
  

\font\tenmsa=msam10
\font\sevenmsa=msam7
\font\fivemsa=msam5
\font\tenmsb=msbm10
\font\sevenmsb=msbm7
\font\fivemsb=msbm5
\newfam\msafam
\newfam\msbfam
\textfont\msafam=\tenmsa  \scriptfont\msafam=\sevenmsa
  \scriptscriptfont\msafam=\fivemsa
\textfont\msbfam=\tenmsb  \scriptfont\msbfam=\sevenmsb
  \scriptscriptfont\msbfam=\fivemsb

\def\hexnumber@#1{\ifnum#1<10 \number#1\else
 \ifnum#1=10 A\else\ifnum#1=11 B\else\ifnum#1=12 C\else
 \ifnum#1=13 D\else\ifnum#1=14 E\else\ifnum#1=15 F\fi\fi\fi\fi\fi\fi\fi}

\def\msa@{\hexnumber@\msafam}
\def\msb@{\hexnumber@\msbfam}
\mathchardef\blacktriangleright="3\msa@49
\mathchardef\blacktriangleleft="3\msa@4A
\catcode`\@=\active

\def\R{\hbox{{\black R}}}
\def\bic{\triangleright\!\!\!\blacktriangleleft}
\def\acti{\triangleleft}

\def\K{{{\cal K}}}
\def\H{{{\cal H}}}
\def\A{{{\cal A}}}
\def\RL{\triangleright\!\!\!\blacktriangleleft}
\def\LR{\blacktriangleright\!\!\!\triangleleft}
\def\LLL{\triangleleft}
\def\RRR{{\bar\triangleright}}
\def\RIMO{\triangleright\!\!\!<}
\def\LECO{>\!\!\!\blacktriangleleft}

\def\balpha{{\bar{\alpha}}}
\def\bbeta{{\bar{\beta}}}
\def\bxi{{\bar{\xi}}}
\def\bpsi{{\bar{\psi}}}
\def\pk{{\cal P}_\kappa}
\def\Upk{{\cal U}_\kappa({\cal P})}
\def\gk{{\cal G}_{\tilde\kappa}}
\def\Ugk{{\cal U}_{\tilde\kappa}({\cal G})}
\def\gmk{\tilde{\cal G}_{(m)\hat\kappa}}
\def\Ugmk{{\cal U}_{\hat\kappa}(\tilde{\cal G}_{(m)})}
\def\Umtr{{\cal U}_{\hat\kappa}(\widetilde{{\cal T}r}_{4})} 
\def\Fgmk{{\rm Fun}_{\hat\kappa}(\tilde G_{(m)})}
\def\Fpk{{\rm Fun}_\kappa(P)}
\def\Fgk{{\rm Fun}_{\tilde\kappa}(G)}
\def\nk{{\cal N}_{\tilde\kappa}}
\def\nmk{\widetilde{\cal N}_{\hat\kappa}}

\def\Fpom{{\rm Fun}_\rho(P(1+1))}
\def\Upom{{\cal U}_\rho({\cal P}(1+1))}
\def\Fgom{{\rm Fun}_\rho(G(1+1))}
\def\Ugom{{\cal U}_\rho({\cal G}(1+1))}
\def\Fgmom{{\rm Fun}_{\hat\rho}(\tilde G_{(m)}(1+1))}
\def\Uph{{\cal U}_h({\cal P})}
\physrev
\nopubblock
\bumpfootsymbolpos

\rightline{FTUV 96--10}
\rightline{IFIC 96--11}
\titlepage
\title{Deformed and extended Galilei group Hopf algebras}
\NAME
\FTUV
\abstract
The $\hat\kappa$-deformed extended Galilei Hopf group algebra,
$\Fgmk$, is introduced. It provides an explicit example of a deformed group 
with cocycle bicrossproduct structure,
and is shown to be the contraction limit of a pseudoextension of the
$\kappa$-Poincar\'e group algebra.
The possibility of obtaining another deformed extended Galilei groups is
discussed, including one obtained from a non-standard Poincar\'e deformation.
\endpage
\chapter{Introduction}
In a recent paper \refmark{\APB} we have discussed the `non-relativistic'
contractions of the deformed $\kappa$-Poincar\'e algebra \refmark{\LNRT},
using the pseudoextension mechanism \refmark{\AACON-\AI} to search for an
{\it extended} deformed Galilei algebra.
In short, this process explains how a trivial (direct product) extension
by the phase group may lead by contraction to a non-trivial central
extension or, in other words, how a two-coboundary may generate a non-trivial 
two-cocycle by contraction.
In the Lie algebra case, the pseudoextension mechanism explains, for
instance, how the direct product ${\cal P}\times u(1)$, where ${\cal P}$
is the Poincar\'e algebra, may lead to the extended Galilei algebra 
$\tilde{\cal G}_{(m)}$; 
other interesting examples may be given both in the
undeformed and the deformed case \refmark{\APBFEST}.

A result of the analysis in \refmark{\APB} is that there are two
possible contractions of the $\kappa$-Poincar\'e algebra \refmark{\LNRT}
$\pk\equiv\Upk$ depending on how the constant $c$ is hidden in $\kappa$,
since the standard $c\to\infty$ limit ($\kappa$ unaltered) leads to the
{\it undeformed} Galilei Hopf algebra ${\cal U}({\cal G})$:

a) $\Ugk$. 
If the deformation parameter $\kappa$ (which has dimensions of inverse
length, $[\kappa]=L^{-1}$) is replaced by $\tilde\kappa/c$ in $\pk$
($[\tilde\kappa]=T^{-1}$), the usual redefinitions ($P_i=X_i\,,\,
P_0\equiv X_t/c \,,\,N_i=cV_i$
with $[X_t]=T^{-1}\,,\,[X]=L^{-1}\,,\,[V]=L^{-1}T$) 
in $\Upk$ lead in the $c\to\infty$
contraction limit to the Hopf algebra $\gk\equiv\Ugk$ given by
$$
\eqalign{
&
[J_i,J_j] =\epsilon_{ijk}J_k\quad,\quad
[J_i,X_j] =\epsilon_{ijk}X_k\quad,\quad
[J_i,X_t] =0\quad, 
\cr 
&
[J_i,V_j] =\epsilon_{ijk}V_k\quad,\quad
[X_i,X_t]=0\quad,\quad
[X_i,X_j]=0\quad,
\cr 
&
[V_i,X_t] =X_i\quad,\quad
[V_i,X_j] =\delta_{ij}{1\over{2\tilde\kappa}}{\bf X}^2-{1\over\tilde\kappa}
X_iX_j\quad,\quad
[V_i,V_j] =0\quad;
\cr
&
\Delta X_t =X_t\otimes 1+1\otimes X_t\quad,\quad
\Delta X_i =X_i\otimes 1 +\exp(-X_t/\tilde\kappa)\otimes X_i\quad,
\cr
&
\Delta J_i =J_i\otimes 1+1\otimes J_i\quad,\quad
\Delta V_i =V_i\otimes 1+ \exp(-X_t/\tilde\kappa)\otimes V_i+{\epsilon_{ijk}
\over\tilde\kappa} X_j\otimes J_k\quad;
\cr
&
S(X_t) =-X_t\quad,\quad
S(X_i) =-\exp(X_t/\tilde\kappa) X_i\quad,
\cr
&
S(J_i) =-J_i \quad,\quad
S(V_i) =-\exp(X_t/\tilde\kappa) V_i +{1\over\tilde\kappa}\epsilon_{ijk}
\exp(X_t/\tilde\kappa) X_jJ_k\quad;
\cr}
\eqn\insi
$$
$\epsilon({\rm all})=0$. This algebra (found in \refmark{\GLODZ} in another
basis) has been shown to have \refmark{\APB} a bicrossproduct
\refmark{\MB} structure: $\Ugk={\cal U}(R\circ B)\bic{\cal
U}_{\tilde\kappa}({\cal T}r_4)=\H\bic\A$ so that 
$\A={\cal U}_{\tilde\kappa}({\cal T}r_4)$ is the
time $(X_t)$ and space $(X_i)$ translations Hopf subalgebra of $\Ugk$ 
and $\H={\cal U}(R\circ B)$ 
is the undeformed Hopf algebra generated by the rotations
$(J_i)$ and boosts $(V_i)$ (here and in the rest of the paper, we shall
use the generic notation $\K=\H\bic\A$ to denote the (right $\RIMO$ / left
$\LECO$) bicrossproduct structure \refmark{\MB} of $\K$).

b) $\Ugmk$.
If $\kappa$ is replaced by $\hat\kappa c$ ($[\hat\kappa]=L^{-2}T$), the
$\kappa$-Poincar\'e algebra $\Upk$ still leads by contraction to the
undeformed Galilei Hopf algebra ${\cal U}({\cal G})$.
However, this redefinition of the deformation parameter allows us to
obtain a non-trivial deformation if the contraction is now performed on a
pseudoextension of $\pk$, $\Upk\times {\cal U}(u(1))$. The result
\refmark{\APB} is the $\hat\kappa$-{\it deformed} extended Galilei algebra
$\gmk\equiv\Ugmk$, where the mass parameter $m$ is introduced, as in the
undeformed case, through the two-coboundary defining the pseudoextension. 
Denoting by $\Xi$ the additional (eleventh) central generator,
the Hopf structure of $\Ugmk$ is given by the
commutators of the undeformed Hopf algebra ${\cal
U}(\tilde{\cal G}_{(m)})$, primitive coproducts etc., but for the exceptions
given below:
$$
\eqalign{
[V_i,X_j]=\delta_{ij}
{\hat\kappa\over 2}(1&-\exp(2m\Xi)/\hat\kappa)\quad;
\cr
\Delta X_i=X_i\otimes 1+\exp(m\Xi/\hat\kappa)\otimes X_i\quad,\quad
&
\Delta V_i=V_i\otimes 1+\exp(m\Xi/\hat\kappa)\otimes V_i\quad;
\cr
S(X_i)=-\exp(-m\Xi/\hat\kappa)X_i\quad,\quad
&
S(V_i)=-\exp(-m\Xi/\hat\kappa)V_i\quad.
\cr}
\eqn\deform
$$
The additional generator $\Xi$ has the dimensions of inverse on an action;
nevertheless it may be rendered dimensionless (in a quantum context) 
multiplying it by the
Planck constant $\hbar$.
The Casimir operators for $\Ugk$ and $\Ugmk$ were also found in
\refmark{\APB}. Since setting $\hat\kappa=0$ in \deform\ the undeformed Hopf
algebra structure of the enveloping algebra ${\cal U}(\tilde{\cal
G}_{(m)})$ is recovered (and, in particular, $[V_i,X_j]=-m\delta_{ij}\Xi$),
we see that the deformation enters in \deform\ only through the central
generator.
In fact, it may be shown \refmark{\APB,\Praga} that $\Ugmk$ has the
(right-left) cocycle bicrossproduct \refmark{\MB} structure
$\Ugmk={\cal U}({\cal G})^\psi\RL_{\xi_{(m)}}{\cal U}(u(1))$, where the right
${\cal U}({\cal G})$-module action ($\RIMO$) $\alpha$ 
(within the generic notation $\H^\psi\RL_\xi\A$, 
$\alpha(a\otimes h)\equiv a\LLL h=a\epsilon_\H(h)$)
of ${\cal G}_{(m)}$
on $u(1)$ is taken to be trivial since $\Xi$ is central in $\Ugmk$
and the left $u(1)$-comodule coaction ($\LECO$) $\beta$
is given by
$$
\beta(X_i)=\exp(m\Xi/\hat\kappa)\otimes X_i\quad,\quad
\beta(V_i)=\exp(m\Xi/\hat\kappa)\otimes V_i
\eqn\dfbeta
$$
($\beta$ is trivial on $X_t$, $J_i\,,$ \ie\ $\beta(X_t)=1\otimes
X_t\,,\,\beta(J_i)=1\otimes J_i$).
The cocycle $\xi_{(m)}$ may be taken to satisfy
\foot{
In \refmark{\APB,\Praga} $\xi_{(m)}$ was taken to be antisymmetric (cf. \dfxi) 
but we wish to be less restrictive here (this freedom is related with the 
election of coboundary as we shall discuss in the next section).}
$$
\xi_{(m)}(V_i,X_j)-\xi_{(m)}(X_j,V_i)= \delta_{ij}{\hat\kappa\over 2}
(1-\exp({2m\Xi\over\hat\kappa}))
\eqn\dfxi
$$
and, finally, $\psi:{\cal U}({\cal G})\to {\cal U}(u(1))\otimes{\cal
U}(u(1))$ is trivial (for $\H^\psi\RL_\xi\A$, triviality
means $\psi(h)=1_\A\otimes 1_\A\epsilon(h)$). 
We note that 
the above structure is not the only one possible.
It may be seen \eg, that
$\Ugmk$ is also the bicrossproduct $\H\bic\A$ of the commutative Hopf
subalgebra $\A=\Umtr$ of $\Ugmk$ generated by 
$<X_t,X_i,\Xi>$ and the undeformed
Hopf algebra $\H={\cal U}(R\circ B)$ 
generated by $<V_i,J_i>$, with right action 
$\alpha:\Umtr\otimes{\cal U}(R\circ B)\to \Umtr$ 
and left coaction 
$\beta:{\cal U}(R\circ B)\to\Umtr\otimes{\cal U}(R\circ B)$
given by
$$
\eqalign{
&
X_t\acti V_i=-X_i\quad,\quad
X_i\acti V_j=-\delta_{ij}{\hat\kappa\over 2}
(1-\exp(2m\Xi/\hat\kappa))\quad,\quad
\Xi\acti V_i=0\quad;
\cr
&
X_t\acti J_i=0\quad,\quad
X_i\acti J_j=\epsilon_{ijk}J_k\quad,\quad
\Xi\acti J_i=0\quad;
\cr
&
\beta(V_i)=\exp(m\Xi/\hat\kappa)\otimes V_i\quad,\quad
\beta(J_i)=1\otimes J_i\quad.
\cr}
\eqn\completar
$$
The existence of more than one structure for a deformed Hopf algebra, as is 
the case for an ordinary Lie algebra, is
not uncommon; see \refmark{\APBFEST}.

Deformed Newtonian and enlarged Newtonian spacetimes may be introduced by
looking at the duals $\nk$ and $\nmk$ of the commutative Hopf subalgebras
${\cal U}_{\tilde\kappa}({\cal T}r_{4})$, $\Umtr$
of $\Ugk$ [\insi] and $\Ugmk$ generated by $<X_t,X_i>$ and
$<X_t,X_i,\Xi>$ respectively.
If a differential calculus covariant under rotations and boosts is now
introduced, it turns out \refmark{\APB} 
that the full commutativity for the
contraction diagrams, which involve the differential calculus
\refmark{\SIT} on the $\kappa$-Minkowski space \refmark{\MR, \LUKRUE, \SZ},
is obtained only for the extended Newtonian spacetime $\nmk$.
Although the addition of a (non-invariant) one-form to
$\Gamma(\nk)$ allows us to solve the rotations/boosts covariance
equations, this new form cannot be interpreted (unlike in $\Gamma(\nmk)$)
as the differential of an additional variable.
This, as the previous discussion, poses the question of whether it is
possible to find an extended deformed Galilei algebra for the
$\tilde\kappa$-deformation as it is the case for $\Ugmk$.
At the same time, it is convenient to have the dual or `group-like'
expressions ($\Fgk$, $\Fgmk$) 
of the above deformed algebras and to look
at the closure of the contraction diagrams involving the deformed
Poincar\'e Hopf group algebra $\Fpk$.
We shall address and answer these questions here, including for
completeness the treatment of the non-standard deformations of the
Poincar\'e group algebra in $(1+1)$ dimensions and its Galilean contractions.
Besides describing possible deformed Galilean groups, our discussion 
will also provide some examples of deformed groups of the cocycle 
bicrossproduct type, a structure of which very few \refmark{\MAJSOB} examples
are known, specially from the group (rather than algebra) point of view.
 
\chapter{The extended deformed Galilei group $\Fgmk$}
Since $\Ugmk$ has a cocycle bicrossproduct \refmark{\MB} structure,
its dual $\Fgmk$ also has this structure. 
In fact, the dual
$K=H_\bxi\LR^\bpsi A$ of $\K=\H^\psi\RL_\xi\A$,
where $H$ and $A$ are the duals of $\H$ and
$\A$, is determined by the dual operations ($\bbeta,\balpha,\bpsi,\bxi$)
of $(\alpha,\beta,\xi,\psi)$.
Quite often the deformation properties of a deformed Hopf algebra $\K$
with the generic structure $\K=\H^\psi\RL_\xi\A$ are mostly described by
those among the mappings 
$\alpha:\A\otimes\H\to\A\,,\,\beta:\H\to\A\otimes\H\,,\,
\xi:\H\otimes\H\to\A\,,\,\psi:\H\to\A\otimes\A$ which 
involve the deformation parameter.
As a result, the dual deformed Hopf group algebra $K=H_\bxi\LR^\bpsi A$
may be found in an easier way by looking for the respective dual maps
($\bbeta:A\to A\otimes H\,,\,\balpha:A\otimes H\to H\, (\balpha:a\otimes h
\mapsto a\RRR h)\,,\,\bpsi:A\to H\otimes H$ and $\bxi:A\otimes A\to H$)
than by direct computation or from `quasiclassical' bialgebra structure
considerations, since $H$ and $A$ are either undeformed or known.
In these cases, the dualization process may be performed in individually
simpler steps.

Since the cocycle bicrossproduct structure of the deformed extended
Galilei algebra $\Ugmk$ was found for $\alpha$ and $\psi$ trivial,
$\bbeta$ and $\bxi$ are also trivial, and only the duals $\balpha$ and
$\bpsi$ of $\beta$ [\dfbeta] and $\xi$ [\dfxi], respectively, need
to be determined since $H={\rm Fun}(G)$ and $A={\rm Fun}(U(1))$ are undeformed.
Since the rotation part in the Galilei group is not important in the
discussion, we shall restrict ourselves initially to 
one spatial dimension.
The commutative, non-cocommutative Hopf algebra structure of 
${\rm Fun}(G(1+1))$
is obvious:
$$
\eqalign{
&
[t,v]=0\quad,\quad
[x,v]=0\quad,\quad
[t,x]=0\quad;
\cr
&
\Delta t=1\otimes t+t\otimes 1\quad,\quad
\Delta x=1\otimes x+x\otimes 1-t\otimes v\quad,\quad
\Delta v=v\otimes +1\otimes v\quad;
\cr
&
S(t)=-t\quad,\quad
S(x)=-x-vt\quad,\quad
S(v)=-v\quad;\quad
\epsilon(t,x,v)=0\quad,\cr}
\eqn\fgal
$$
where the group algebra generators $t$ ($x$) 
correspond to the time (space) translation and $v$ to the boost. 
The duality relations are
$$
<X,x>=1\quad,\quad
<X_t,t>=1\quad,\quad
<V,v>=1\quad,
\eqn\duality
$$
the others being zero.
To determine $\Fgmk$ we introduce an additional generator $\chi$ dual to
$\Xi$, $<\Xi,\chi>=1$, with $\Delta\chi=\chi\otimes
1+1\otimes\chi\,,\,S(\chi)=-\chi\,,\,\epsilon(\chi)=0$. 
Then, using eqs. \dfbeta\ and relations
such as
$<V,\chi\RRR v>\equiv <V,\balpha(\chi\otimes v)>=<\beta(V),\chi\otimes v>$,
the left action $\balpha$ ($\RRR$) dual to $\beta$ is immediately seen to be
$$
\chi\RRR t=0\quad,\quad
\chi\RRR x={m\over\hat\kappa}x\quad,\quad
\chi\RRR v={m\over\hat\kappa}v\quad.
\eqn\dfrac
$$

Finding the dual cocycle $\bpsi$ of $\xi$ requires a more careful
analysis. Since, as mentioned, 
explicit examples \refmark{\MAJSOB,\APB,\APBFEST} of cocycle 
bicrossproduct constructions for deformed Hopf algebras are rather scarce,
we shall provide it with some detail.
Recalling that 
$\xi:{\cal U}({\cal G})\otimes{\cal U}({\cal G})\to{\cal U}(u(1))$,
let us take
$$
\xi(V,X)=B(\Xi)\quad,\quad
\xi(X,V)=B(\Xi)-{\hat\kappa\over 2}(1-\exp({2m\Xi\over\hat\kappa}))\quad
\eqn\dfb
$$
(so that \dfxi\ is satisfied), where $B(\Xi)$ is as yet undetermined.
Since $[V,X_t]=X$, it follows that
$$
\eqalign{
&
\xi(V,VX_t)-\xi(V,X_tV)=B(\Xi)\quad,
\cr
&
\xi(VX_t,V)-\xi(X_tV,V)=B(\Xi)-
{\hat\kappa\over 2}(1-\exp({2m\Xi\over\hat\kappa}))
\quad.
\cr}
\eqn\dfbii
$$
We may now use the two-cocycle condition 
\refmark{\DOI,\MB} on $\xi$ to obtain more information.
Since the right action $\alpha$ ($\LLL$) was trivial, this condition
reduces to
$$
\xi(h_{(1)}g_{(1)}\otimes f)\xi(h_{(2)}\otimes g_{(2)})=
\xi(h\otimes g_{(1)}f_{(1)})\xi(g_{(2)}\otimes f_{(2)})
\eqn\twococi
$$
where $h,g,f\in{\cal U}({\cal G})$.
Assuming the natural normalization
$\xi(h\otimes 1)=\xi(1\otimes h)=1\epsilon(h)$
and $\xi(1\otimes 1)=1$, the application of $\twococi$ to the Galilei
algebra generators, which have primitive coproducts, reduces it to
$\xi(hg\otimes f)=\xi(h\otimes gf)$. Thus, we have in particular
$$
\xi(VX_t,V)=\xi(V,X_tV)\equiv A(\Xi)\quad,
\eqn\dfa
$$
introducing a new unknown function.
Then, with
$$
<B(\Xi),\chi>=B\quad,\quad
<A(\Xi),\chi>=A\quad,
\eqn\dfab
$$
we find, besides eqs. \dfbii, \dfa, the relations
$$
\eqalign{
&
\xi(V,VX_t)=A(\Xi)+B(\Xi)=\xi(V^2,X_t)\quad,
\cr
&
\xi(X_t,V^2)=A(\Xi)-B(\Xi)+
{\hat\kappa\over 2}(1-\exp({2m\Xi\over\hat\kappa}))\quad.
\cr}
\eqn\dfabii
$$
Now, using that 
$<{\hat\kappa\over 2}(1-\exp({2m\Xi\over\hat\kappa})),\chi>=-m\,,\,
<VX_t,vt>=1\,,\, <V^2,v^2>=2,$ the dual $\bpsi$ of $\xi$ is found to be
$$
\eqalign{
\bpsi(\chi)=
&
A[{1\over 2}v^2\otimes t+vt\otimes v+v\otimes vt+{1\over 2}t\otimes v^2]+
\cr
&
B[v\otimes x+x\otimes v+v\otimes vt+{1\over 2}v^2\otimes t-{1\over
2}t\otimes v^2]+
\cr
&
m[x\otimes v-{1\over 2} t\otimes v^2]\quad.
\cr}
\eqn\dualcoci
$$
The last bracket is recognized as a form of the Galilei non-trivial
two-cocycle \refmark{\LL}; 
since the first two depend on the constants $A$, $B$, they
must correspond to two-coboundaries.
This may be trivially checked in a Lie group language if we think of the 
terms in the r.h.s. of \dualcoci\ as the product of the unprimed and primed
group parameters defining the standard group two-cocycle, $\omega(g,g')$
say. 
In this way, it is seen that the first two brackets are generated by the
one-cochains $\eta(g)=-{1\over 2}v^2t$ and $\eta(g)=-{1\over 2}v^2t-vx$,
respectively, trough
$\omega_{\rm cob}(g,g')=\eta(g)+\eta(g')-\eta(gg')$.
Thus, we do not need worrying here about the cocycle condition for $\bpsi$, 
although we shall come back to it in sec. 3.
We shall only note now that, in the present Hopf algebra context, the
condition expressing that $\bpsi$ is a coboundary 
is translated into
$$
\eqalign{
\bpsi_{\rm cob}(\chi)=
&
1\otimes\gamma(\chi)+\gamma(\chi^{(1)})\otimes\chi^{(2)}-\Delta(\gamma(\chi))
\cr
=
&
1\otimes\gamma(\chi)+\gamma(\chi)\otimes 1-\Delta(\gamma(\chi))\quad,
\cr}
\eqn\coboun
$$
where $\gamma$ is a linear mapping $\gamma: A \to H$ which is convolution
invertible (\ie, there exists $\gamma^{-1}$ such that 
$\gamma(a_{(1)})\gamma^{-1}(a_{(2)})=\gamma^{-1}(a_{(1)})\gamma(a_{(2)})=
\epsilon(a)\ $);
$\gamma(\chi)$ is given by $-{1\over 2}v^2t\,$ [$-{1\over 2}v^2t-vx$] for the
first [second] bracket in \dualcoci.
Ignoring the coboundaries, we thus find that $\Fgmk$ is defined by eqs.
\fgal\ to which one has to add those dictated by \dfrac\ \ie
$$
[\chi,t]=0\quad,\quad
[\chi,x]={m\over\hat\kappa}x\quad,\quad
[\chi,v]={m\over\hat\kappa}v\quad,
\eqn\dfcomm
$$
the two cocycle \dualcoci, the last term of which 
modifies the primitive coproduct
$\Delta\chi$ to read
$$
\Delta\chi=\chi\otimes 1+1\otimes\chi+m(x\otimes v-{1\over 2}t\otimes v^2)
\quad,
\eqn\dfcop
$$
and the antipode and counit of $\chi$,
$$
S(\chi)=-\chi+m(xv+{1\over 2}tv^2)\quad,\quad
\epsilon(\chi)=0\quad.
\eqn\dfant
$$
As expected from $\Ugmk$, the non-commutative nature of $\Fgmk$ only shows
up in the commutation properties \dfcomm\ of the additional generator $\chi$.

Moving now to $(1+3)$ dimensions by including the rotations is not difficult.
The only modified coproducts and antipodes are
$$
\eqalign{
&
\Delta x_i=1\otimes x_i+x_j\otimes R^j_{~i}-t\otimes v_i\quad,\quad
\Delta v_i=1\otimes v_i+v_j\otimes R^j_{~i}\quad;
\cr
&
S(x_i)=-x_j(R^{-1})^j_{~i}-tv_j(R^{-1})^j_{~i}\quad,\quad
S(v_i)=-v_j(R^{-1})^j_{~i}\quad,
\cr}
\eqn\inrot
$$
and equations such as \duality, \dfcomm, \dfcop\ and \dfant\ acquire 
the appropriate vector indices
(for instance, $S(\chi)=-\chi+m[{\bf xv}+{1\over 2}t{\bf v}^2]$.)
As for $R^j_{~i}$, it commutes with all other group algebra generators,
$(\Delta R)^i_{~j}=R^i_{~k}\otimes R^k_{~j}$,
and its presence modifies the cocycle to read
$$
\bpsi(\chi)=m(x_i\otimes R^i_{~j}v^j-{1\over 2}t\otimes {\bf v}^2)\quad.
\eqn\dualrotcoci
$$
If the angle-like generator $\varphi=\chi/\hbar$ is used, $\Fgmk$ may be
rightly called the {\it deformed quantum-mechanical Galilei group}.

\chapter{$\Fgmk$ as a contraction of a pseudoextension of $P_{\hat\kappa}$ 
(${\rm Fun}_{\hat\kappa}(P)$)}
Consider the $\kappa$-Poincar\'e group algebra $P_\kappa\equiv\Fpk$ in
$(1+1)$ dimensions \refmark{\ZAUGG, \SZ}. It is defined by the relations
$$
[x_1,x_0]={x_1\over\kappa}\quad,\quad
[\alpha,x_0]={1\over\kappa}\sinh\alpha\quad,\quad
[\alpha,x_1]=-{1\over\kappa}(\cosh\alpha-1)\quad,
\eqn\poicom
$$
(where $x_0=ct$ ($x_1$) refers to time (space) and $\alpha$ characterizes the
boost in the $x$ direction) and by
$$
\eqalign{
&
\Delta\alpha=\alpha\otimes 1+1\otimes\alpha\quad,\quad
\Delta x_0=1\otimes x_0+x_0\otimes\cosh\alpha-x_1\otimes\sinh\alpha\quad,
\cr
&
\Delta x_1=1\otimes x_1-x_0\otimes\sinh\alpha+x_1\otimes\cosh\alpha\quad;
\cr
&
S(\alpha)=-\alpha\quad,\quad
S(x_0)=-x_0\cosh\alpha-x_1\sinh\alpha\quad,\quad
\cr
&
S(x_1)=-x_1\cosh\alpha-x_0\sinh\alpha\quad;\quad
\epsilon(\alpha,x_0,x_1)=0\quad.
\cr}
\eqn\dualpoi
$$
Introduce now the commutative, co-commutative Hopf algebra generated by 
a new generator $\hat\chi$ with dimensions of an action, such that
$$
[\hat\chi,{\rm all}]=0\quad;\quad
\Delta\hat\chi=\hat\chi\otimes1+1\otimes\hat\chi\quad,\quad
S(\hat\chi)=-\hat\chi\quad;\quad
\epsilon(\hat\chi)=0\quad.
\eqn\dfchia
$$
If we now make the redefinition $\chi=\hat\chi-mcx_0$ using the
one-cochain $mcx_0$ (which diverges in the contraction limit) it follows that
$$
\eqalign{
&
\Delta\chi=\chi\otimes 1+1\otimes\chi-
mc[x_0\otimes(\cosh\alpha-1)-x_1\otimes\sinh\alpha]\quad,
\cr
&
S(\chi)=-\chi+mc[x_0(\cosh\alpha-1)+x_1\sinh\alpha]\quad.
\cr}
\eqn\dfchi
$$
In terms of $\chi$, the commutators in \dfchia\ imply
$$
[\chi,x_0]=0\quad,\quad
[\chi,x_1]={mc\over\kappa}x_1\quad,\quad
[\chi,\alpha]={mc\over\kappa}\sinh\alpha\quad,
\eqn\chicomm
$$
and we see again that in order to have a sensible contraction limit we
need substituting $\hat\kappa c$ for $\kappa$, in which case eqs. \dfcomm\
and \dfcop\ are recovered ($\alpha\sim v/c$) and hence $\Fgmk$.
Thus, there is full duality and closure among the algebra-like
contraction diagrams for $\Ugmk$ \refmark{\APB} and the group-like ones
for $\Fgmk$.

A question that naturally arises is whether we can construct a similar
extended group 
algebra from the Galilei deformation $\Fgk$ dual of $\Ugk$, eqs.
\insi. As $\Ugk$, $\Fgk$ has a bicrossproduct structure so that the
dualization (in two dimensions) of ${\cal U}(B)\RL{\cal
U}_{\tilde\kappa}({\cal T}r_2)$ leads immediately to
${\rm Fun}_{\tilde\kappa}(G(1+1))$ as (cf. \fgal)
$$
\eqalign{
&
[t,x]=-{1\over\tilde\kappa}x\quad,\quad
[x,v]={v^2\over 2\tilde\kappa}\quad,\quad
[t,v]=-{v\over\tilde\kappa}\quad;
\cr
&
\Delta t=t\otimes 1+1\otimes t\quad,\quad
\Delta x=x\otimes 1+1\otimes x-t\otimes v\quad,\quad
\Delta v=v\otimes 1+1\otimes v\quad;
\cr
&
S(t,x,v)=(-t,-x-vt,-v)\quad,\quad
\epsilon(t,x,v)=0\quad.
\cr}
\eqn\fkgal
$$
Let us add a ${\rm Fun}(U(1))$ factor to ${\rm Fun}_{\tilde\kappa}(G)$ (a
similar discussion could be presented for $\Ugk$).
Since we want the corresponding algebra generator to be central in $\Ugk$,
$\alpha$ has to be trivial and hence its dual $\bbeta$ is also trivial.
The two-cocycle condition for the map
$\bpsi:A\to H\otimes H$ (here $A={\rm Fun}(U(1))\,,\, 
H={\rm Fun}_{\tilde\kappa}(G)$) may be found as a consequence of the 
coassociativity requirement. It is given by \refmark{\MB}
$$
\Delta\bpsi(a_{(1)})^{(1)}\bpsi(a_{(2)}^{~(1)})\otimes\bpsi(a_{(1)})^{(2)}
a_{(2)}^{~(2)}=
\bpsi(a_{(1)})^{(1)}\otimes\Delta\bpsi(a_{(1)})^{(2)}\bpsi(a_{(2)})\quad,
\eqn\dualcocicond
$$
where the lower indices refer to the coproduct and the upper ones refer
either to
$\bbeta$ ($\bbeta(a)=a^{(1)}\otimes a^{(2)}\in A\otimes H$) or to the
components of $\bpsi$ ($\bpsi(a)=\bpsi(a)^{(1)}\otimes\bpsi(a)^{(2)})$.
Since the right coaction $\bbeta$ is trivial in our case, 
$a^{(1)}=a\,,\,a^{(2)}=1$.
Moreover, since $\bpsi(1)=1_H\otimes 1_H$ and 
$\Delta a=a\otimes 1+1\otimes a\,$ \ie,
$\Delta\chi=\chi\otimes 1+1\otimes\chi$,
eq. \dualcocicond\ reduces to
$$
\Delta\bpsi(\chi)^{(1)}\otimes\bpsi(\chi)^{(2)}+\bpsi(\chi)\otimes 1=
\bpsi(\chi)^{(1)}\otimes\Delta\bpsi(\chi)^{(2)}+1\otimes\bpsi(\chi)
\eqn\dualcocycle
$$
(it is easy to check, we note in passing, that \dualrotcoci\ for $\Fgmk$
satisfies condition \dualcocycle).
Thus, in the search for a $\Fgk$ we have to look for a $\bpsi(\chi)$ which
in the undeformed limit must reduce to the last term in \dualcoci.
The usual dimension assignments (which we have consistently kept) indicate 
that a $\tilde\kappa$-deformation 
of the Galilei two-cocycle (last bracket in \dualcoci) may be
described by
$$
\bpsi_{\tilde\kappa}=\bpsi(\chi)(1\otimes f({mv^2\over\hbar\tilde\kappa}))
\quad.
\eqn\dfdeforcoci
$$
Omitting the constants we may now impose the two-cocycle condition
\dualcocycle\ to \dfdeforcoci\ written as
$\bpsi_{\tilde\kappa}\propto x\otimes vf(v^2)-{1\over 2}t\otimes v^2f(v^2)$.
This leads to
$$
\eqalign{
&
\Delta x\otimes v f(v^2)-{1\over 2}\Delta t\otimes v^2f(v^2)+x\otimes vf(v^2)
\otimes 1-{1\over 2}t\otimes v^2f(v^2)\otimes 1=
\cr
&
x\otimes\Delta(vf(v^2))-{1\over 2} t\otimes \Delta(v^2f(v^2))
+1\otimes x\otimes vf(v^2)
-{1\over 2}1\otimes t\otimes v^2f(v^2)\quad,
\cr}
\eqn\demosi
$$
which implies
$$
\eqalign{
&
\Delta(vf(v^2))=1\otimes vf(v^2)+vf(v^2)\otimes 1\quad,
\cr
&
\Delta(v^2f(v^2))=1\otimes v^2f(v^2)+v^2f(v^2)\otimes 1+2v\otimes vf(v)\quad.
\cr}
\eqn\demosii
$$
These equations are inconsistent with the form of $\Delta v$ in \fkgal\ {\it 
and} among each other unless $f$ is constant.
Thus, $\bpsi$ is unmodified by $\tilde\kappa$ and $\Delta\chi$
is again given by \dfcop.
We may now try to complete the commutation relations in \fkgal\ with those
for $\chi$ by imposing that the coproduct (as given by \fkgal, \dfcop) is
an algebra homomorphism. 
Notice that the addition of a two-cocycle does not modify the 
expressions \fkgal; 
only the commutators involving $\chi$ and $\Delta\chi$ need to be found.
The result is that there is no solution in the 
presence of $\tilde\kappa$ if the 
constant $f$ is non-zero
\foot{
This result disagrees with that in \refmark{\GGKM} inasmuch as
the projective representations of $\Fgk$ are related to a two-cocycle 
Hopf group algebra extension (which is not discussed there).}.
This agrees with the fact that no 
$\tilde\kappa$-deformation of the Hopf algebra ${\cal U}(\tilde{\cal G}_{(m)})$
can be obtained from $\Upk\times{\cal U}(u(1))$ by contraction
\refmark{\APB}.

\chapter{Structure of the non-standard (1+1) deformed Poincar\'e group and
their Galilean contractions}

A non-standard deformation $\Uph$ of the Poincar\'e Hopf algebra 
may be obtained by contraction from the
${\cal U}_h(sl(2,\R))$ deformation \refmark{\Ohn,\AKS} and has been recently
studied \refmark{\KSAA,\BCHOS,\BHOPS}.
We shall show first that it has a bicrossproduct structure and then study
its Galilean contractions.

In a `light-cone' basis $\Upom$ ($\rho$ is the parameter remaining after
the contraction limit $\epsilon\to 0$ which is performed after setting
$h=\epsilon\rho$, $[\epsilon]=L^{-1}=[\rho^{-1}]$) 
may be written in the form
$$
\eqalign{
&
[N,P_+]={1-\exp(-2\rho P_+)\over 2\rho}\quad,\quad
[N,P_-]=-P_-\quad,\quad
[P_+,P_-]=0\quad;
\cr
&
\Delta P_+=P_+\otimes 1+1\otimes P_+\quad,\quad
\Delta P_-=\exp(-2\rho P_+)\otimes P_-+P_-\otimes 1\quad,
\cr
&
\Delta P_-=\exp(-2\rho P_+)\otimes N+N\otimes 1\quad;\quad
S(N)=-\exp(2\rho P_+)N\quad,
\cr
&
S(P_+)=-P_+\quad,\quad
S(P_-)=-\exp(2\rho P_+)P_-\quad;\quad\epsilon(N,P_{\pm})=0\quad.
\cr}
\eqn\npoal
$$
$\Upom$ has a bicrossproduct structure $\K=\H\RL\A$ where in
this case $\H$ is generated by $N$ with primitive coproduct and $\A$ is
the Abelian Hopf subalgebra of $\Upom$ generated by the translations
$P_\pm$, the right action $\alpha:\A\otimes\H\to\A$
and left coaction $\beta:\H\to\A\otimes\H$ being given, respectively, by
$$
P_+\LLL N=-{1-\exp(-2\rho P_+)\over 2\rho}\quad,\quad
P_-\LLL N=P_-\quad;\quad
\beta(N)=\exp(-2\rho P_+)\otimes N\quad.
\eqn\npobic
$$
We may easily see \eg, that the formulae which give the coproduct and
antipode in $\K$ \refmark{\MB},
$$
\eqalign{&
\Delta_\K(h\otimes a)=h_{(1)}\otimes h_{(2)}^{~(1)}a_{(1)}
\otimes h_{(2)}^{~(2)}\otimes a_{(2)}\quad;
\cr
&
S(h\otimes a)=(1_\H\otimes S_\A(h^{(1)}a))(S_\H(h^{(2)})\otimes 1_{\A})\quad,
\cr}
\eqn\copant
$$
immediately reproduce $\Delta N$ and $S(N)$ ($N$ is represented in $\K$ by
$N\otimes 1$ and, in eq. \copant, $N^{(1)}=\exp(-2\rho P_+)\,,\,N^{(2)}=N$
by \npobic\ etc.)
It is simple to construct by duality ($<P_{\pm},x_{\pm}>=1$) the associated
non-standard $(1+1)$ spacetime Hopf algebra, which is defined by the relations
$$
[x_+,x_-]=-2\rho x_-\quad;\quad
\Delta x_\pm=x_\pm\otimes 1+1\otimes x_\pm\quad;\quad
S(x_\pm)=-x_\pm\;;\;
\epsilon(x_\pm)=0\;.
\eqn\npospace
$$

We are interested here, however, in constructing the whole dual Hopf algebra
$\Fgom$ and its possible extension, and in obtaining them from $\Fpom$.
Denoting the variable dual to $N$ by $\alpha$, $<N,\alpha>=1$, and using
the fact that $\Upom$ is a bicrossproduct, the duals $\bbeta$ of $\alpha$
and $\balpha\, (\RRR)$ of $\beta$ are found to be
$$
\bbeta(x_\pm)=x_\pm\otimes e^{\mp\alpha}\quad;\quad
x_+\RRR\alpha=-2\rho(1-e^{-\alpha})\quad,\quad
x_-\RRR\alpha=0\quad,
\eqn\dfaccion
$$
(to find \eg, $x_+\RRR\alpha$ one needs considering $\beta$ for powers of $N$,
$\beta(N^m)$). In this way the {\it non-standard Poincar\'e Hopf algebra}
$\Fpom$ is found to be
$$
\eqalign{
&
[x_+,\alpha]=-2\rho(1-e^{-\alpha})\quad,\quad
[x_-,\alpha]=0\quad,\quad
[x_+,x_-]=-2\rho x_-\quad;
\cr
&
\Delta x_+=1\otimes x_++x_+\otimes e^{-\alpha}\quad,\quad
\Delta x_-=1\otimes x_-+x_-\otimes e^\alpha\quad,
\cr
&
\Delta \alpha=\alpha\otimes 1+1\otimes\alpha\quad;\quad
S(x_+)=-x_+e^\alpha\quad,
\cr &
S(x_-)=-x_-e^{-\alpha}\quad,\quad
S(\alpha)=-\alpha\quad;\quad
\epsilon(x_\pm,\alpha)=0\quad,
\cr}
\eqn\npogroup
$$
which reproduces \refmark{\BHOPS,\KSAA}.

We now find the {\it non-standard deformed Galilei group} $\Fgom$. In
terms of the standard 
($x_0=x_++x_-\,,\,x_1=x_+-x_-$) basis,
$\Fpom$ is given by
$$
\eqalign{
&
[x_0,\alpha]=-2\rho(1-e^{-\alpha})\quad,\quad
[x_1,\alpha]=-2\rho(1-e^{-\alpha})\quad,
\cr
&
[x_0,x_1]=2\rho(x_0-x_1)\quad;
\cr
&
\Delta\alpha=\alpha\otimes 1+1\otimes\alpha\quad,\quad
\Delta x_0=1\otimes x_0+x_0\otimes\cosh\alpha-x_1\otimes\sinh\alpha\quad,
\cr
&
\Delta x_1=1\otimes x_1-x_0\otimes\sinh\alpha+x_1\otimes\cosh\alpha\quad;
\cr
&
S(\alpha)=-\alpha\quad,\quad
S(x_0)=-x_0\cosh\alpha-x_1\sinh\alpha\quad,
\cr
&
S(x_1)=-x_1\cosh\alpha-x_0\sinh\alpha\quad;\quad
\epsilon(\alpha,x_0,x_1)=0\quad,
\cr}
\eqn\npogroupx
$$
which reproduces, in the coalgebra sector, the standard $P(1+1)$ group law.
The non-standard deformation $\Fgom$ now follows from the $c\to\infty$ limit 
for $x_0=ct\,,\,x_1=x$ and $\alpha\sim v/c$, with
the result
$$
\eqalign{
&
\Delta t=1\otimes t+t\otimes 1\quad,\quad
\Delta x=1\otimes x+x\otimes 1-t\otimes v\quad,\quad
\Delta v=1\otimes v+v\otimes 1\quad;
\cr
&
[t,v]=0\quad,\quad
[x,v]=-2\rho v\quad,\quad
[t,x]=2\rho t\quad;
\cr
&
S(t)=-t\quad,\quad
S(x)=-x-vt\quad,\quad
S(v)=-v\quad.
\cr}
\eqn\ngalgroup
$$
The dual algebra, $\Ugom$ is easily found to be
$$
\eqalign{
&
[X_t,X]=0\quad,\quad
[X_t,V]=-{1\over4\rho}(1-\exp(-4\rho X))\quad,\quad
[X,V]=0\quad;
\cr
&
\Delta X_t=X_t\otimes 1+\exp(-2\rho X)\otimes X_t\quad,\quad
\Delta X=1\otimes X+X\otimes 1\quad,
\cr
&
\Delta V=V\otimes 1+\exp(-2\rho X)\otimes V\quad;\quad
S(X_t)=-\exp(2\rho X)X_t\quad,
\cr
&
S(X)=-X\quad,\quad
S(V)=-\exp(2\rho X)V\quad;\quad
\epsilon(X_t,X,V)=0\quad;
\cr}
\eqn\ngalalg
$$
It is not a new Hopf algebra;
it is the deformed Heisenberg-Weyl algebra ${\cal U}_\rho(HW)$
(the quantum Heisenberg group $H_q(1)$ of \refmark{\CGSTa}), and it has both a
bicrossproduct and a cocycle bicrossproduct structure \refmark{\APBFEST}.

To complete the picture, it is interesting to close the contraction
diagrams by obtaining $\Ugom\ (={\cal U}_\rho(HW))$, eqs. \ngalalg, by
contracting $\Upom$ in \npoal\ by means of the standard redefinitions.
It turns out, however, that the na\"{\i}ve change of basis 
($P_0={1\over 2}(P_++P_-)\,,\,P_1={1\over 2}(P_+-P_-)$)
from light cone to
standard variables is not adequate, and that a more complicated one 
(which may be justified \refmark{\BHOPS} in terms of $T$ matrix \refmark{\FG}
considerations) is
required to perform the contraction, namely, $P_+=P_0+P_1\,,\, P_-={1\over
2\rho}[\exp(-2\rho(P_0+P_1))+1-2\exp(-2\rho P_0)]$.
In terms of these $P_{0,1}$ generators the algebra $\Upom$ is more
complicated than in \npoal, but may be seen to lead to $\Ugom$, eqs.
\ngalalg, in the contraction limit.

Finally, we may look for a {\it non-standard extended Galilei group}.
This can be obtained following the by now familiar procedure, which
involves the addition of $\hat\chi$ as in \dfchia, and the redefinition
$\hat\chi=\chi+mcx_0$. 
Then \npogroupx\ leads to \dfchi\ and to
$$
[\chi,x_0]=0\quad,\quad
[\chi,x_1]=-2mc\rho(x_0-x_1)\quad,\quad
[\chi,\alpha]=2mc\rho(1-e^{-\alpha})\quad.
\eqn\dfcommut
$$
Eqs. \dfcommut\ show that the contraction requires a redefinition
of the deformation constant, $\rho=\hat\rho/c^2\,,\,[\hat\rho]=L^3T^{-2}$.
This leads to $\Fgmom$ defined by eqs. \fgal, \dfcop\ together
with 
$$
[\chi,t]=0\quad,\quad
[\chi,x]=-2m\hat\rho t\quad,\quad
[\chi,v]=0\quad.
\eqn\dfcommutator
$$
$\Fgmom$ is a very mild deformation, only manifest in the $[\chi,x]$ 
commutator.

\ack
This paper has been partially supported by the CICYT research grant AEN93-187.
One of us (JCPB) wishes to acknowledge a FPI grant from the Spanish Ministry 
of Education and Science and the CSIC.

\refout
\end